\documentclass[12pt, a4paper, oneside]{amsart}
\usepackage{amssymb, amsmath, amsthm, verbatim, amsbsy,cite}
\usepackage[english]{babel}
\usepackage{amsaddr}
\newtheorem{theorem}{Theorem}[section]
\newtheorem{conj}[theorem]{Conjecture}
\newtheorem{problem}[theorem]{Problem}

\newcommand{\inv}{\operatorname{inv}}
\newcommand{\Aut}{\operatorname{Aut}}

\newcommand{\cd}{\operatorname{cd}}

\newcommand{\cs}{\operatorname{cs}}
\newcommand{\Cl}{\operatorname{Cl}}
\newcommand{\eo}{\operatorname{eo}}

\newcommand{\mC}{\mathbb{C}}

\begin{document}

\title{Arithmetic invariants for finite simple and related groups}

\author{Andrey V. Vasil'ev}
\address{Sobolev Institute of Mathematics, Novosibirsk, Russia}
\email{vasand@math.nsc.ru}

\begin{abstract} In this short note, we address problems concerning the characterization of simple and related groups by various arithmetic invariants. We propose a uniform approach to such questions and discuss both known results and open problems in this area.
\smallskip

\smallskip

\end{abstract}

\maketitle

\section{Introduction}

One may view finite group theory as the study of the interplay between the algebraic structure of a group and its arithmetic properties, that is, properties expressed in terms of numerical parameters. The well-known theorems of Lagrange, Sylow, and Feit--Thompson are among numerous examples.

Here we are interested in the case where a set of arithmetic parameters of a finite group $G$ gives rise to a system of invariants that allow us to distinguish $G$ from all other finite groups up to isomorphism. We focus our attention on the following three sets:
\begin{itemize}
\item $\eo(G)=\{|x| : x\in G\}$ the set of element orders;
\item $\cs(G)=\{|C| : C\in\Cl(G)\}$ the set of conjugacy class sizes;
\item $\cd(G)=\{\chi(1): \chi\in \operatorname{Irr}(G)\}$ the set of complex irreducible character degrees.
\end{itemize}

We denote by $\eo^*(G), \cs^*(G), \cd^*(G)$  the corresponding multisets, with multiplicities taken into account.
\smallskip

As examples show, one cannot expect $\eo(G), \cs(G), \cd(G)$ or even $\eo^*(G), \cs^*(G), \cd^*(G)$ to provide a full system of invariants in the class of all finite groups. However, the situation changes for finite simple groups and groups related to them, such as almost simple and quasisimple groups.

The importance of simple groups in finite group theory is the same as the importance of prime numbers in number theory. The problem of characterizing them within the class of all finite groups has become an active area of research since the classification theorem was announced and the complete list of finite simple groups became available. Recall that, according to the classification of finite simple groups (CFSG), the finite simple groups are precisely
\begin{enumerate}
\item the groups of prime order;
\item the alternating groups of degree at least $5$;
\item the simple classical groups;
\item the simple exceptional groups of Lie type;
\item the 26 sporadic groups.
\end{enumerate}

We begin by discussing three well-known conjectures. The first two arose in 1987 in a private communication between Shi and Thompson (see, e.g., \cite{25Shi}), while the third was posed by Huppert in 2000, see \cite{Huppert2000}. All three of them were later added to the {\em Kourovka Notebook} \cite[Problems 12.38, 12.39, 18.99]{KT}.

\begin{conj}{\em(Shi's Conjecture)}\label{c:Shi}
Let $L$ be a simple group and let $G$ be a group such that $\eo(G) = \eo(L)$ and $|G|=|L|$. Then $G\simeq L$.
\end{conj}

\begin{conj}{\em(Thompson's Conjecture)}\label{c:Thompson}
Let $L$ be a nonabelian simple group and let $G$ be a group such that $\cs(G) = \cs(L)$ and $Z(G)=1$. Then $G\simeq L$.
\end{conj}

\begin{conj}{\em(Huppert's Conjecture)}\label{c:Hupppert}
Let $L$ be a nonabelian simple group and let $G$ be a group such that $\cd(G) = \cd(L)$. Then $G\simeq L \times  A$, where $A$ is abelian.
\end{conj}

\textbf{Remark.} It is clear that the validity of Huppert’s conjecture would imply that $G\simeq L$ if one adds either of the two conditions $|G|=|L|$ or $Z(G)=1$ from the Shi and Thompson conjectures to the condition $\cd(G) = \cd(L)$. Moreover, the validity of Thompson's conjecture would imply that $G\simeq L$ if the condition $Z(G)=1$ in Thompson’s conjecture were replaced by $|G|=|L|$ (see, \cite[Lemma~2.4]{GV2026}). It is also clear that, if we assume that $|G|=|L|$, then $L$ need not be assumed to be nonabelian in any of the three statements.
\smallskip

Shi's conjecture was proved in \cite{VGM09}. Thompson’s conjecture was proved for all finite simple groups except the alternating groups in~\cite{19GorS}; for the alternating groups, it was proved modulo the binary Goldbach conjecture in~\cite{19GorA}. Huppert’s conjecture has been verified for all finite simple groups except the classical groups; see \cite{25TV} for references. Similar problems are also worth considering for groups related to simple groups.

The purpose of this note is to provide a uniform framework for thinking about questions of this kind. Recall that an {\em invariant} of an equivalence relation on a set $M$ is a function
$$f: M\rightarrow Q,$$
where $Q$ is some set (usually, a field or the ring of integers), such that $f$ is constant on each equivalence class. Since $f(x)\neq f(y)$ implies that $x,y\in M$ are non-equivalent, invariants can be used to distinguish equivalence classes. A system of invariants $\mathcal{F}$ is {\em full} if, for every non-equivalent $x,y\in M$, there exists $f\in\mathcal{F}$ such that $f(x)\neq f(y)$.

Given a subset $N\subseteq M$ closed under the equivalence relation, we say that a system of invariants $\mathcal{F}$ is {\em full for $N$}, if for every $z\in N$ and an arbitrary $x\in M$ non-equivalent to~$z$, there is $f\in\mathcal{F}$ such that $f(x)\neq f(z)$.

In our case, the equivalence is the isomorphism relation on the set of finite groups and we are looking for systems of invariants full for simple groups and for groups related to them.

\section{Invariants for simple groups}

First, we fix $n=n(G)=|G|$, the order of a group $G$, which is certainly an invariant of the isomorphism relation on the set of finite groups.

Let $D(G)=D(n)$ be the set of the divisors of~$n$ and $\tau(n)=|D(n)|$. It is well known that for $\inv\in\{\eo,\cs,\cd\}$, $\inv(G)\subseteq D(G)=D(n)$. Arrange the elements $d_i$, $i=1,\ldots,\tau(n),$ of $D(n)$ in ascending order. For $\inv\in\{\eo,\cs,\cd\}$ and every $d_i\in D(G)$, set
$$
\delta^{\inv}_i(G):=
\begin{cases}
1, \text{ if } d_i\in\inv(G) \\
0, \text{ otherwise}.
\end{cases}
$$
For $\inv^*\in\{\eo^*,\cs^*,\cd^*\}$, and every $d_i\in D(G)$, set
$$
\delta^{\inv^*}_i(G):=\text{multiplicity of }d_i\text{ in }\inv^*(G).
$$
It is clear that the functions $\delta^{\inv}_i(G)$ and $\delta^{\inv^*}_i(G)$ are invariants of $G$ for every $i=1,\ldots,\tau(n)$.
Let $\delta^{\inv}(G)=(\ldots\delta^{\inv}_{i}(G)\ldots)$ and $\delta^{\inv^*}(G)=(\ldots\delta^{\inv^*}_{i}(G)\ldots)$  be the tuples of length $\tau(n)$.\smallskip

In this language, Shi’s conjecture can be reformulated as follows.

\begin{theorem}\label{t:Shi}
The system of invariants $(n(G),\delta^{\eo}(G))$ is full for the simple groups.
\end{theorem}

Recently in \cite{GV2026}, we proved

\begin{theorem}\label{t:new}
The system of invariants $(n(G),\delta^{\cs}(G))$ is full for the alternating and symmetric groups.
\end{theorem}

Thompson’s conjecture for all simple groups other than the alternating groups, together with \cite[Lemma~2.4]{GV2026}, implies that this system of invariants is full for all simple groups (see also \cite[Corollary~1.3]{GV2026}).

\begin{theorem}\label{t:ThompsonOrder}
The system of invariants $(n(G),\delta^{\cs}(G))$ is full for the simple groups.
\end{theorem}

In view of these results, it is natural to ask the following question.

\begin{problem}\label{prob:ThompsonOrder}
Is the system of invariants $(n(G),\delta^{\cd}(G))$ full for the simple groups?
\end{problem}

\textbf{Remark.} It is clear that Huppert’s conjecture implies a positive answer to Problem~\ref{prob:ThompsonOrder}, but not conversely. Since Huppert's conjecture has been verified for all simple groups except the classical groups of dimension greater than~4, it suffices to solve the problem for the latter.
\smallskip

Since
$$
n(G)=\sum_{d_i\in D(G)}\delta^{\eo^*}_i(G)=\sum_{d_i\in D(G)}d_i\,\delta^{\cs^*}_i(G)=\sum_{d_i\in D(G)}d_i^2\,\delta^{\cd^*}_i(G),
$$
the multiset $\inv^*(G)$ determines the systems of invariants $(n(G),\delta^{\inv^*}(G))$ and $(n(G),\delta^{\inv}(G))$. Although it is not known whether the system of invariants $(n(G),\delta^{\cd}(G))$ is full for all simple groups, it follows from \cite{10TV, 12TV1, 12TV} that $(n(G),\delta^{\cd^*}(G))$ is full. Thus, the following general assertion on the characterization of simple groups holds true.

\begin{theorem}\label{t:MultSimple}
For $\inv^*\in\{\eo^*,\cs^*,\cd^*\}$, the system of invariants $(n(G),\delta^{\inv^*}(G))$ is full for the simple groups.
\end{theorem}

\section{Groups related to simple groups}

In many cases, one needs to deal with groups that are not simple but are closely related to simple groups. We recall that
\begin{itemize}
\item $G$ is \emph{almost simple} if $S\leq G\leq\Aut S$ for a nonabelian simple group $S$ (e.g., $G=S_m$ with $S=A_m$),
\item $G$ is \emph{quasisimple} if $G=[G,G]$ and $G/Z(G)=S$ for a nonabelian simple group $S$ (e.g., $G=\operatorname{SL}_m(q)$ with $S=\operatorname{PSL}_m(q)$),
\item $G$ is \emph{almost quasisimple} if $G$ has a quasisimple normal subgroup $H$ with $C_G(H)\leq H$.
\end{itemize}
In what follows, in each of these cases we say that $G$ is \emph{related} to $S$ (in the last case, $S=H/Z(H)$).\smallskip

The following general problem arises.

\begin{problem}\label{p:general}
For which groups $G$ related to a finite simple group, and for which $\inv\in\{\eo,\cs,\cd\}$, is the system of invariants $(n(G),\delta^{\inv}(G))$ or $(n(G),\delta^{\inv^*}(G))$ full in the class of all finite groups{\em?}
\end{problem}

We begin our discussion of this problem with the case of symmetric groups. We know from Theorem~\ref{t:new} that the system of invariants $(n(G),\delta^{\cs}(G))$ is full for them. The same holds true for $\inv=\eo$.

\begin{theorem}{\em\cite[Theorems~6,7]{90Bi}}\label{t:Bi}
The system of invariants $(n(G),\delta^{\eo}(G))$ is full for the symmetric groups.
\end{theorem}

\textbf{Remark.} It is worth noting that the main result of~\cite{90Bi} (available only in Chinese) can now be deduced from a much stronger assertion. Namely, it follows from \cite{14Gor,16GorGr} that, for every positive integer $m\not\in\{2, 3, 4, 5, 6, 8, 10\}$, the symmetric group $S_m$ is uniquely characterized by the set $\eo(S_m)$ in the class of all finite groups. In the exceptional cases, one can verify that the equality $(n(G),\delta^{\eo}(G))=(n(S_m),\delta^{\eo}(S_m))$ yields $G\simeq S_m$.
\smallskip

The case $\inv=\cd$ remains open.

\begin{problem}\label{prob:CdSym}
Is the system of invariants $(n(G),\delta^{\cd}(G))$  full for the symmetric groups?
\end{problem}

Nevertheless, it follows from the main result of~\cite{11TV} that, if we replace $(n(G),\delta^{\cd}(G))$ by $(n(G),\delta^{\cd^*}(G))$, we obtain a system of invariants full for the symmetric groups. Therefore, the following holds.

\begin{theorem}\label{t:sym}
For $\inv^*\in\{\eo^*,\cs^*,\cd^*\}$, the system of invariants $(n(G),\delta^{\inv^*}(G))$ is full for the symmetric groups.
\end{theorem}

In what follows, it is more convenient to discuss the three cases $\inv=\eo,$ $\cs$, and $\cd$ separately. We begin with the case $\inv=\cd$.

The irreducible character degrees of a group $G$, counted with multiplicities, form the first column of the character table of~$G$ and uniquely determine the complex group algebra~${\mC}G$. This makes the characterization of groups related to simple groups by $\cd^*(G)$ an active area of research. In our terminology, the following holds true, see~\cite{15BNOTV} and references therein.

\begin{theorem}\label{t:quasicd*}
The system of invariants $(n(G),\delta^{\cd^*}(G))$ is full for all quasisimple groups.
\end{theorem}

In the case of almost simple groups, only partial results are available. Besides the symmetric groups mentioned above, it is known that the system of invariants $(n(G),\delta^{\cd^*}(G))$ is   full for almost simple groups with socle $PSL_2(q)$ \cite{17HA}, for the groups $PGU_3(q)$ \cite{17SI}, and for groups $G$ with $PSL_n(q)<G\leq PGL_n(q)$ \cite{22SI}, or $PSU_n(q)<G\leq PGU_n(q)$, where $q+1$ divides neither $n$ nor $n-1$ \cite{20SIS}.

Turning to the case where irreducible character degrees are considered without multiplicities, we note that there have been many attempts to formulate analogues of Huppert’s conjecture for groups related to simple groups, see, e.g., \cite{15HMTVW,20SI,23Al} and the partial results obtained there. Within our framework, we suggest unifying them by adding $n(G)$ as an additional invariant. Although it is not known whether the system $(n(G),\delta^{\cd}(G))$ is full for all simple groups, even the following question remains open.

\begin{problem}\label{prob:CdRel}
Does there exist a quasisimple group $G$ for which the system of invariants $(n(G),\delta^{\cd(G)})$ is not full?
\end{problem}

Very little is known about the characterization of groups related to simple groups by conjugacy class sizes. Apart from Theorem~\ref{t:new} on symmetric groups, it is known (see~\cite{22Gor}) that $\cs(G)=\cs(\operatorname{SL}_2(5))$ implies that $G\simeq\operatorname{SL}_2(5)\times A$ for an abelian group $A$. In particular, the system of invariants $(n(G),\delta^{\cs}(G))$ is full for $\operatorname{SL}_2(5)$.
As in the previous case, the following question remains open.

\begin{problem}\label{prob:CSRel}
Does there exist a quasisimple group $G$ for which the system of invariants $(n(G),\delta^{\cs(G)})$ is not full?
\end{problem}

We now turn to the case $\inv=\eo$. Although the system of invariants $(n(G),\delta^{\eo}(G))$ is full for the symmetric groups $S_m$, this is not true for all almost simple and quasisimple groups related to $A_m$. Recall that $A_6\simeq\operatorname{PSL}_2(9)$, $|\Aut A_6|=4\cdot|A_6|$, and as one can easily verify using the {\em Atlas of Finite Groups} \cite{Atlas}:
$$
\eo(\Aut A_6)=\eo(C_2\times\operatorname{PGL}_2(9))=\eo(C_2\times\operatorname{SL}_2(9)).
$$

\textbf{Remark.} Among almost simple groups with socle $A_m$, the case $m=6$ is the only exceptional one, since $\Aut A_m=S_m$ for $m\neq2,6$.

There are infinitely many groups related to simple classical groups for which the system of invariants $(n(G),\delta^{\eo}(G))$ is not full. For example, according to  \cite[Proposition~2.6]{08But}, for groups related to $\operatorname{PSL}_m(q)$, where $m$ is odd and $m\neq p^\ell+1$, we have
$$\eo(\operatorname{PGL}_m(p^k))=\eo(\operatorname{SL}_m(p^k))\mbox{ and }|\operatorname{PGL}_m(p^k)|=|\operatorname{SL}_m(p^k)|.$$
Similar examples exist for other classical groups.

Moreover, even the system of invariants $(n(G),\delta^{\eo^*}(G))$ is not full for all almost simple groups. As was observed by Thompson (see, e.g., \cite[Section~4.3]{25Shi}), there are two maximal subgroups $G_1\simeq L_3(4).2$ and $G_2\simeq 2^4\rtimes A_7$ of the Mathieu group $M_{23}$ such that $\eo^*(G_1)=\eo^*(G_2)$.

Despite these examples, we believe that $n(G)$ and $\eo(G)$  are “almost sufficient” to characterize most groups related to simple groups. We propose the following conjecture.

\begin{conj}\label{c:involutions}
The order $n(G)$ of $G$, the number $i(G)$ of involutions in $G$, and $\delta^{\eo}(G)$ form a system of invariants  that is full for all but finitely many almost quasisimple groups.
\end{conj}

As noted in the Introduction, even $\eo^*(G)$, $\cs^*(G)$, $\cd^*(G)$ cannot provide full systems of invariants within the class of all finite groups. Nevertheless, it is natural to ask what restrictions the equality $\inv^*(G)=\inv^*(H)$, or weaker conditions such as $\inv(G)=\inv(H)$, imposes on the internal structures of groups $G$ and $H$?

For example, we know that the equality $\eo^*(G)=\eo^*(H)$ guarantees neither that both groups $G$ and~$H$ are solvable nor that they are both semisimple. Indeed, on the one hand, Piwek recently constructed in \cite{Piw24} a solvable group $G$ and nonsolvable group $H$ with $\eo^*(G)=\eo^*(H)$, see also \cite{Mul24} for smaller examples, thus solving a long-standing question of Thompson, see \cite[Problem~12.37]{KT}. On the other hand, the subgroup $G_1\simeq L_3(4).2$ of the Mathieu group $M_{23}$  is semisimple, i.e., $G_1$ has trivial solvable radical, whereas the subgroup $G_2\simeq 2^4\rtimes A_7$ of the same group obviously does not, though $\eo^*(G_1)=\eo^*(G_2)$. We wonder whether there are similar examples for $\inv=\cs$ or $\inv=\cd$.

\begin{problem}\label{prob:solv}
Let $\inv\in\{\cs,\cd\}$, and let $\inv^*(G)=\inv^*(H)$.
\begin{enumerate}
\item If $G$ is solvable, must $H$ be solvable?
\item If $G$ has trivial solvable radical, must $H$ have trivial solvable radical?
\end{enumerate}
\end{problem}

\section{Concluding Remarks}

Results on arithmetic characterizations of simple and related groups have applications in algebra and representation theory. As noted above, $\cd^*(G)$ uniquely determines the complex group algebra ${\mC}G$. Thus, the fact that $\cd^*(G)$ provides a set of invariants full for all simple and even for all quasisimple groups (Theorem~\ref{t:quasicd*} above) yields an affirmative answer to Question~2 in Brauer's seminal paper \cite{63Br} for quasisimple groups, see also the introduction in \cite{15BNOTV}. In view of \cite[Corollary 5.2]{08KimLuRC}, Shi's conjecture (Theorem~\ref{t:Shi}) implies  that a finite simple group and a finite group with the same Burnside rings are isomorphic, thereby answering Yoshida's question in \cite[Problem 2]{87Yos} for finite simple groups. Let us also mention that Theorem~\ref{t:Shi} has recently appeared in research on the {\em group isomorphism problem} in computational complexity~\cite{BS2022,GL25} and even in machine learning~\cite{HJCS25}.

It is worth noting that, in all these applications, the order of the group is included among the input invariants. To some extent, this justifies our approach of regarding $(n(G),\delta^{\inv}(G))$ as a {\em minimal} set of invariants for a group~$G$, while regarding $(n(G),\delta^{\inv^*}(G))$ as a {\em maximal} one.

\medskip

\textbf{Acknowledgments.} This article grew out of a talk delivered during my visit to the Sino-Russian Mathematics Center at Peking University.
The author is deeply grateful to Professor Zhang Jiping, Director of the SRMC, for the invitation, warm hospitality, and fruitful discussions. The research was carried out within the framework of the Sobolev Institute of Mathematics project FWNF-2026-0017.

\end{document}